% From shelah@math.wisc.edu  Fri Apr 19 20:53:52 1996
% Received: from conley.math.wisc.edu (conley.math.wisc.edu [144.92.166.10]) by sunset.ma.huji.ac.il (8.6.11/8.6.10) with SMTP id UAA23013 for <shlhetal@sunrise.huji.ac.il>; Fri, 19 Apr 1996 20:53:03 +0300
% Received: by conley.math.wisc.edu;
%           id AA15400; 4.1/42; Fri, 19 Apr 96 12:18:37 CDT
% Date: Fri, 19 Apr 96 12:18:37 CDT
% From: Saharon Shelah <shelah@math.wisc.edu>
% Message-Id: <9604191718.AA15400@conley.math.wisc.edu>
% To: shlhetal@sunrise.ma.huji.ac.il
% Subject: 484 (with lie kecheng
% X-Status: 
% Status: RO
% 
\ifx\shlhetal\undefinedcontrolsequence\let\shlhetal\relax\fi

\magnification=1200
\overfullrule0pt
\input mssymb
\expandafter\ifx\csname bib4plain.tex\endcsname\relax
  \expandafter\gdef\csname bib4plain.tex\endcsname{}
\else \message{Hey!  Apparently you were trying to \string twice.   This does not make sense.}
\errmessage{Please edit your file (probably \jobname.tex) and remove
any duplicate ``\string\input'' lines} \fi

%  This file should be inputted if you want to use 
%  bibtex fom within plain TeX. 
      % Not really need for standard
       % bibtex files, but these commands
\def\renewcommand{\newcommand}	       % are used in our literal-unsrt.bst
\edef\cite{\the\catcode`@}%
\catcode`@ = 11
\let\@oldatcatcode = \cite
\chardef\@letter = 11
\chardef\@other = 12
%
%
% Next come some things that will be useful later.
%
% Make an outer definition into an inner one (due to Chris Thompson).
% The arguments should be the control sequence to be defined, and the
% new of the \outer control sequence, as characters; the control
% sequence #1 is defined to be just the same as \csname#2\endcsname, but
% not \outer.  For example, \@innerdef\innernewcount{newcount} would
% define \innernewcount to be a non-outer version of \newcount.
%
\def\@innerdef#1#2{\edef#1{\expandafter\noexpand\csname #2\endcsname}}%
%
% We use \@innerdef to make some of our allocations local, because
% Eplain includes our code inside a conditional.  We put @'s in the
% names to minimize the (already small) chance of conflicts.
%
\@innerdef\@innernewcount{newcount}%
\@innerdef\@innernewdimen{newdimen}%
\@innerdef\@innernewif{newif}%
\@innerdef\@innernewwrite{newwrite}%
%
%
% Swallow one parameter.
%
\def\@gobble#1{}%
%
%
% Use TeX 3.0's \inputlineno to get the line number, for better error
% messages, but if we're using an old version of TeX, don't do anything.
%
\ifx\inputlineno\@undefined
   \let\@linenumber = \empty % Pre-3.0.
\else
   \def\@linenumber{\the\inputlineno:\space}%
\fi
%
%
% The following macro \@futurenonspacelet (from the TeXbook) behaves
% essentially like \futurelet except that it discards any implicit or
% explicit space tokens that intervene before a nonspace is scanned:
%
\def\@futurenonspacelet#1{\def\cs{#1}%
   \afterassignment\@stepone\let\@nexttoken=
}%
\begingroup % The grouping here avoids stepping on an outside use of `\\'.
\def\\{\global\let\@stoken= }%
\\ % now \@stoken is a space token (\\ is a control symbol, so that
   % space after it is seen).
\endgroup
\def\@stepone{\expandafter\futurelet\cs\@steptwo}%
\def\@steptwo{\expandafter\ifx\cs\@stoken\let\@@next=\@stepthree
   \else\let\@@next=\@nexttoken\fi \@@next}%
\def\@stepthree{\afterassignment\@stepone\let\@@next= }%
%
%
% \@getoptionalarg\CS gets an optional argument from the input, enclosed
% in brackets, then expands \CS.  We set \@optionalarg to \empty if we
% don't find one, otherwise to the text of the argument.  This assumes
% the brackets don't have a funny category code.
%
\def\@getoptionalarg#1{%
   \let\@optionaltemp = #1%
   \let\@optionalnext = \relax
   \@futurenonspacelet\@optionalnext\@bracketcheck
}%
%
% The \expandafter's in this macro let us avoid the use of \aftergroup,
% which is somewhat more expensive.
%
\def\@bracketcheck{%
   \ifx [\@optionalnext
      \expandafter\@@getoptionalarg
   \else
      \let\@optionalarg = \empty
      % We can't do the \temp after the \fi, because then the \temp gets
      % in the way of reading the optional argument from the input, if
      % we do have one.
      \expandafter\@optionaltemp
   \fi
}%
\def\@@getoptionalarg[#1]{%
   \def\@optionalarg{#1}%
   \@optionaltemp
}%
%
%
% From LaTeX.
%
\def\@nnil{\@nil}%
\def\@fornoop#1\@@#2#3{}%
\def\@for#1:=#2\do#3{%
   \edef\@fortmp{#2}%
   \ifx\@fortmp\empty \else
      \expandafter\@forloop#2,\@nil,\@nil\@@#1{#3}%
   \fi
}%
\def\@forloop#1,#2,#3\@@#4#5{\def#4{#1}\ifx #4\@nnil \else
       #5\def#4{#2}\ifx #4\@nnil \else#5\@iforloop #3\@@#4{#5}\fi\fi
}%
\def\@iforloop#1,#2\@@#3#4{\def#3{#1}\ifx #3\@nnil
       \let\@nextwhile=\@fornoop \else
      #4\relax\let\@nextwhile=\@iforloop\fi\@nextwhile#2\@@#3{#4}%
}%
%
%
% This macro tests if a file \jobname.#1 exists, and sets \if@fileexists
% appropriately.  If an optional argument is given, it is used as the
% root part of the filename instead of \jobname.
%
\@innernewif\if@fileexists
\def\@testfileexistence{\@getoptionalarg\@finishtestfileexistence}%
\def\@finishtestfileexistence#1{%
   \begingroup
      \def\extension{#1}%
      \immediate\openin0 =
         \ifx\@optionalarg\empty\jobname\else\@optionalarg\fi
         \ifx\extension\empty \else .#1\fi
         \space
      \ifeof 0
         \global\@fileexistsfalse
      \else
         \global\@fileexiststrue
      \fi
      \immediate\closein0
   \endgroup
}%
%
%
%% [[[start of BibTeX-specific stuff]]]
%
% Now come the four main LaTeX commands and their associated .aux
% commands.  Just as in LaTeX, \bibliographystyle defines the BibTeX
% style name (.bst file, that is), and \bibliography defines the
% database (.bib) file(s).  The corresponding .aux-file commands are
% \bibstyle and \bibdata, which are there only for BibTeX's (but not
% LaTeX's) use.
%
\def\bibliographystyle#1{%
   \@readauxfile
   \@writeaux{\string\bibstyle{#1}}%
}%
\let\bibstyle = \@gobble
%
% As well as writing the \bibdata command to tell BibTeX which .bib
% files to read, we read the .bbl file that BibTeX (or a person,
% conceivably) has produced.  We use \bblfilebasename as the root of the
% filename to read; this defaults to \jobname.
%
\let\bblfilebasename = \jobname
\def\bibliography#1{%
   \@readauxfile
   \@writeaux{\string\bibdata{#1}}%
   \@testfileexistence[\bblfilebasename]{bbl}%
   \if@fileexists
      % We just output a non-discardable item (the `whatsit' with the
      % \bibdata command).  This means that the glue that will be
      % inserted next (\parskip or \baselineskip, most likely) will be a
      % legal breakpoint.  Most likely, this is after some kind of
      % heading, where we don't want to allow a page break.  So:
      \nobreak
      \@readbblfile
   \fi
}%
\let\bibdata = \@gobble
%
% The \nocite{label,label,...} command writes its argument to \@auxfile,
% unless instructed not to, but produces no text in the document.  Both
% the \nocite and \cite commands produce \citation commands in the .aux file.
%
\def\nocite#1{%
   \@readauxfile
   \@writeaux{\string\citation{#1}}%
}%
\@innernewif\if@notfirstcitation
%
% \cite[note]{label,label,...} produces the citations for the labels as
% well.  If the optional argument `note' is present, it's added after
% the labels.  Since \cite calls \nocite to do its .aux-file writing,
% \cite doesn't need to call \@readauxfile (\nocite does).
%
\def\cite{\@getoptionalarg\@cite}%
%
% Typeset the citations for the labels in #1, followed by the note, if
% it exists.  To change the citation's format in the text, redefine one
% or more `\print...' macros, whose defaults appear later in this file.
%
\def\@cite#1{%
   % Remember the optional argument, in case one of the macros we call
   % below ends up looking for an optional argument itself.  For
   % example, if a \cite[note] triggers reading the .aux file, then the
   % [note] would be clobbered, since \@testfileexistence looks for an
   % optional arg.
   \let\@citenotetext = \@optionalarg
   % Start printing the text, beginning with a left bracket by default.
   \printcitestart
   % It's complicated, but because \nocite puts a `whatsit' onto the list,
   % \nocite should follow \printcitestart.  It's conceivable, but very
   % unlikely, that this `whatsit' will cause a problem (glue that doesn't
   % disappear when you want it to is the most likely symptom), requiring
   % a change either to \printcitestart or to the label that the .bst file
   % produces.
   \nocite{#1}%
   \@notfirstcitationfalse
   \@for \@citation :=#1\do
   {%
      \expandafter\@onecitation\@citation\@@
   }%
   \ifx\empty\@citenotetext\else
      \printcitenote{\@citenotetext}%
   \fi
   \printcitefinish
}%
\def\@onecitation#1\@@{%
   \if@notfirstcitation
      \printbetweencitations
   \fi
   \expandafter \ifx \csname\@citelabel{#1}\endcsname \relax
      \if@citewarning
         \message{\@linenumber Undefined citation `#1'.}%
      \fi
      % Give it a dummy definition:
      \expandafter\gdef\csname\@citelabel{#1}\endcsname{%
% Change: marginal remark added, goldstrn@math.huji.ac.il, 
% goldstern@tuwien.ac.at, May 1996 mg
%  !!! change !!!
\strut
\vadjust{\vskip-\dp\strutbox
\vbox to 0pt{\vss\parindent0cm \leftskip=\hsize 
\advance\leftskip3mm
\advance\hsize 4cm\strut\openup-4pt 
\rightskip 0cm plus 1cm minus 0.5cm ?  #1 ?\strut}}
         {\tt
            \escapechar = -1
            \nobreak\hskip0pt
            \expandafter\string\csname#1\endcsname
            \nobreak\hskip0pt
         }%
      }%
   \fi
   % Now produce the text, whether it was undefined or not.
   \csname\@citelabel{#1}\endcsname
   \@notfirstcitationtrue
}%
%
% Given a label `foo', the macro `\b@foo' is supposed to
% hold the text that should be produced.
%
\def\@citelabel#1{b@#1}%
%
% So, how does a citation label get defined?  When we read the .bbl file
% (below), a \bibitem writes out a \@citedef command.  And when we read
% the \@citedef, we define \@citelabel{#1}, where #1 is the user's
% label.
%
\def\@citedef#1#2{\expandafter\gdef\csname\@citelabel{#1}\endcsname{#2}}%
%
%
% Reading the .bbl file also produces the typeset bibliography.  Please
% notice, however, that we do not produce the title for the references
% (e.g., `References'), as LaTeX does.  The formatting and spacing of
% that title, whether it should go into the headline, and so on, are all
% things determined by your format.  We cannot know those things in
% advance.  If you wish, you can define \bblhook to produce the title.
% Or just do it before the \bibliography command.
%
\def\@readbblfile{%
   % Define a counter to tell us which item number we are on, unless
   % we've already defined it (because the document has more than one
   % bibliography).
   \ifx\@itemnum\@undefined
      \@innernewcount\@itemnum
   \fi
   \begingroup
      \def\begin##1##2{%
         % ##1 is just `thebibliography'.
         % ##2 is the widest label.
         % We set (new dimen) \biblabelwidth based on the widest label
         \setbox0 = \hbox{\biblabelcontents{##2}}%
         \biblabelwidth = \wd0
      }%
      \def\end##1{}% ##1 is `thebibliography' again.
      %
      % Here we have two possibilities:
      % \bibitem[typesetlabel]{citationlabel}
      % \bibitem{citationlabel}
      % If we have the second of these, the citations are numbered, starting
      % from one; we use our own count register \@itemnum for this.
      %
      \@itemnum = 0
      \def\bibitem{\@getoptionalarg\@bibitem}%
      \def\@bibitem{%
         \ifx\@optionalarg\empty
            \expandafter\@numberedbibitem
         \else
            \expandafter\@alphabibitem
         \fi
      }%
      \def\@alphabibitem##1{%
         % Need \xdef here for various reasons.
         \expandafter \xdef\csname\@citelabel{##1}\endcsname {\@optionalarg}%
         % Left-justify alpha labels, unless \biblabel{pre,post}contents
         % are already defined.
         \ifx\biblabelprecontents\@undefined
            \let\biblabelprecontents = \relax
         \fi
         \ifx\biblabelpostcontents\@undefined
            \let\biblabelpostcontents = \hss
         \fi
         \@finishbibitem{##1}%
      }%
      \def\@numberedbibitem##1{%
         \advance\@itemnum by 1
         \expandafter \xdef\csname\@citelabel{##1}\endcsname{\number\@itemnum}%
         % Right-justify numeric labels, unless \biblabel{pre,post}contents
         % are already defined.
         \ifx\biblabelprecontents\@undefined
            \let\biblabelprecontents = \hss
         \fi
         \ifx\biblabelpostcontents\@undefined
            \let\biblabelpostcontents = \relax
         \fi
         \@finishbibitem{##1}%
      }%
      \def\@finishbibitem##1{%
         \biblabelprint{\csname\@citelabel{##1}\endcsname}%
         \@writeaux{\string\@citedef{##1}{\csname\@citelabel{##1}\endcsname}}%
         \ignorespaces
      }%
      %
      % Do the printing (we're producing the bibliography, remember).
      %
      \let\em = \bblem
      \let\newblock = \bblnewblock
      \let\sc = \bblsc
      % Punctuation won't affect spacing;
      \frenchspacing
      % the penalties below are from LaTeX's [article,book,report].sty;
      \clubpenalty = 4000 \widowpenalty = 4000
      % the next two values come from LaTeX's \sloppy command;
      \tolerance = 10000 \hfuzz = .5pt
      \everypar = {\hangindent = \biblabelwidth
                      \advance\hangindent by \biblabelextraspace}%
      \bblrm
      % the \parskip is a guess at what looks good;
      \parskip = 1.5ex plus .5ex minus .5ex
      % and the space between label and text comes from LaTeX's \labelsep.
      \biblabelextraspace = .5em
      \bblhook
      \input \bblfilebasename.bbl
   \endgroup
}%
%
% The widest label's width is useful for redefining \biblabelprint;
% you redefine \biblabelwidth, in effect, by redefining the
% \biblabelcontents macro that appears below.  And \biblabelextraspace,
% which is redefinable inside \bblhook, is added to \biblabelwidth to
% determine the amount of hanging indentation.
%
\@innernewdimen\biblabelwidth
\@innernewdimen\biblabelextraspace
%
% Now come the main macros that are related to the printing of the
% bibliography.  Since you might want to redefine them, they are given
% default definitions outside of \@readbblfile.
%
% The first one controls the printing of a bibliography entry's label.
% If you change it, make sure that it starts with something like
% \noindent or \indent or \leavevmode that puts TeX into horizontal mode
% (even if the label itself is empty); otherwise, the hanging
% indentation will get messed up in certain circumstances.
%
\def\biblabelprint#1{%
   \noindent
   \hbox to \biblabelwidth{%
      \biblabelprecontents
      \biblabelcontents{#1}%
      \biblabelpostcontents
   }%
   \kern\biblabelextraspace
}%
%
% If you are using numeric labels, and you want them left-justified
% (numeric labels by default are right-justified), do something like:
%     \def\biblabelprecontents{\relax}
%     \def\biblabelpostcontents{\hss}
%
% By default the labels are typeset in \bblrm, and enclosed in brackets.
%
\def\biblabelcontents#1{{\bblrm [#1]}}%
%
% The main text, too, is typeset using \bblrm, which is \rm by default.
%
\def\bblrm{\rm}%
%
% Emphasis for producing, e.g., titles, is done with \it by default.
%
\def\bblem{\it}%
%
% Some styles use a caps-and-small-caps font for author names.  LaTeX
% defines an \sc command but plain TeX doesn't, so we need one here.
% The definition below doesn't load the font unless it's needed, but it
% tries to load only the 10pt version, because it might not exist at
% other point sizes.
%
\def\bblsc{\ifx\@scfont\@undefined
              \font\@scfont = cmcsc10
           \fi
           \@scfont
}%
%
% The major parts of an entry are separated with \bblnewblock.  The
% numbers below are taken from LaTeX's `article' style.
%
\def\bblnewblock{\hskip .11em plus .33em minus .07em }%
%
% Here's where you stick any other bibliography-formatting goodies, or
% redefine the values above.
%
\let\bblhook = \empty
%
%
% Here are the four default definitions for formatting the in-text
% citations.  These are what you redefine (after your \input btxmac but
% before your \bibliography) to get parens instead of brackets, or
% superscripts, or footnotes, or whatever.
%
\def\printcitestart{[}%         left bracket
\def\printcitefinish{]}%        right bracket
\def\printbetweencitations{, }% comma, space
\def\printcitenote#1{, #1}%     comma, space, note (if it exists)
%
% That scheme is pretty flexible.  For example you could use
%     \def\printcitestart{\unskip $^\bgroup}
%     \def\printcitefinish{\egroup$}
%     \def\printbetweencitations{,}
%     \def\printcitenote#1{\hbox{\sevenrm\space (#1)}}
%     \font\eighttt = cmtt8
%     \scriptfont\ttfam = \eighttt
% to get superscripted in-text citations.  (The scriptfont stuff
% exists only to print an undefined citation; it's in cmtt8 because
% there is no cmtt7.)  To get something radically different, however,
% you'll have to define your own \cite command.
%
% When we read `\citation' from the .aux file, it means nothing.
%
\let\citation = \@gobble
%
%
% Now comes the stuff for dealing with LaTeX's \newcommand.  As
% mentioned earlier, this \newcommand will redefine a preexisting
% command; that's different from how LaTeX's \newcommand behaves.
%
\@innernewcount\@numparams
%
% \newcommand{\foo}[n]{text} defines the control sequence \foo to have
% n parameters, and replacement text `text'.
%
\def\newcommand#1{%
   \def\@commandname{#1}%
   \@getoptionalarg\@continuenewcommand
}%
%
% Figure out if this definition has parameters.
%
\def\@continuenewcommand{%
   % If no optional argument, we have zero parameters.  Otherwise, we
   % have that many.
   \@numparams = \ifx\@optionalarg\empty 0\else\@optionalarg \fi \relax
   \@newcommand
}%
%
% \@numparams is how many arguments this command has.  The name of the
% command is \@commandname.  The replacement text for the new macro is #1.
%
\def\@newcommand#1{%
   \def\@startdef{\expandafter\edef\@commandname}%
   \ifnum\@numparams=0
      \let\@paramdef = \empty
   \else
      \ifnum\@numparams>9
         \errmessage{\the\@numparams\space is too many parameters}%
      \else
         \ifnum\@numparams<0
            \errmessage{\the\@numparams\space is too few parameters}%
         \else
            \edef\@paramdef{%
               % This is disgusting, but \loop doesn't work inside \edef,
               % because \body isn't defined.
               \ifcase\@numparams
                  \empty  No arguments.
               \or ####1%
               \or ####1####2%
               \or ####1####2####3%
               \or ####1####2####3####4%
               \or ####1####2####3####4####5%
               \or ####1####2####3####4####5####6%
               \or ####1####2####3####4####5####6####7%
               \or ####1####2####3####4####5####6####7####8%
               \or ####1####2####3####4####5####6####7####8####9%
               \fi
            }%
         \fi
      \fi
   \fi
   \expandafter\@startdef\@paramdef{#1}%
}%
%
%% [[[end of BibTeX-specific stuff]]]
%
%
% Names of references (arguments given in the \cite and \nocite
% commands) and file names (arguments given in the \bibliography and
% \bibliographystyle commands) are recorded in \jobname.aux, called the
% \@auxfile in these macros.  Here's how they get read in.
%
\def\@readauxfile{%
   \if@auxfiledone \else % remember: \@auxfiledonetrue if \noauxfile is defined
      \global\@auxfiledonetrue
      \@testfileexistence{aux}%
      \if@fileexists
         \begingroup
            % Because we might be in horizontal mode when \@readauxfile
            % is called (if it's in response to a \cite or \nocite), we
            % want to ignore all the would-be spaces at the ends of
            % lines in the aux file.  Fortunately, it's highly unlikely
            % an end-of-line might actually be desired.
            % And because we don't change the category code of anything
            % but @, primitives like \gdef can't be used to define labels
            % in the aux file.  The solution adopted by btxmac.tex is to
            % write `\@citedef{LABEL}{DEFINITION}' to the aux file, and
            % use \csname on LABEL.
            \endlinechar = -1
            \catcode`@ = 11
            \input \jobname.aux
         \endgroup
      \else
         \message{\@undefinedmessage}%
         \global\@citewarningfalse
      \fi
      \immediate\openout\@auxfile = \jobname.aux
   \fi
}%
%
% The \@readauxfile macro does all that work the first time it's called.
% Since it's called once for every \cite, \nocite, \bibliography, and
% \bibliographystyle command that the user issues, we need to remember
% whether the work's been done.  It's considered done if we're not to do
% it---that is, if \noauxfile is defined.
%
\newif\if@auxfiledone
\ifx\noauxfile\@undefined \else \@auxfiledonetrue\fi
%
% It's conceivable you'd want to change how other characters are read;
% to do that, change their category code before doing \input btxmac.
%
%
% After reading the .aux file, \@readauxfile opens it for writing.
% The \@writeaux macro does the actual writing (as long as
% \noauxfile is undefined).
%
\@innernewwrite\@auxfile
\def\@writeaux#1{\ifx\noauxfile\@undefined \write\@auxfile{#1}\fi}%
%
%
% A macro package that uses btxmac.tex might define
% \@undefinedmessage (before doing an \input btxmac).
%
\ifx\@undefinedmessage\@undefined
   \def\@undefinedmessage{No .aux file; I won't give you warnings about
                          undefined citations.}%
\fi
%
% Even if citations are undefined, we want to complain only if
% \@citewarningtrue.  The default is to set \@citewarningtrue unless
% \noauxfile is defined.  Again, a macro package that uses
% btxmac.tex might want to redefine this.
%
\@innernewif\if@citewarning
\ifx\noauxfile\@undefined \@citewarningtrue\fi
%
%
% Finally, before leaving we restore @'s old category code.
%
\catcode`@ = \@oldatcatcode

\hsize=16.5truecm
\vsize=23truecm

\def\HH{{\cal H}}

\def\FF{{\cal F}}

\def\ra{{\rangle}}
\def\la{{\langle}}
\def\alo{{\aleph_\omega}}
\def\n{{\noindent}}
 2
\font\tit=cmr10 scaled \magstep2
\font\mtit=cmsy10 scaled \magstep2
\font\bit=cmr10 scaled \magstep1

\def\n{{\noindent}}

\def\mod{{\rm ~mod~}}
\def\Suc{{\rm Suc}}

\def\Sup{{\rm Sup}}
\def\Min{{\rm Min}}
\def\Max{{\rm Max}}
\def\II{{\cal I}}
\def\TT{{\cal T}}
\def\PP{{\cal P}}
\def\FF{{\cal F}}
\def\HH{{\cal H}}
\def\mo{{\displaystyle {\mathop {=}_{\rm def}}}}

\vglue1.25cm
\centerline{\tit Cofinalities of Elementary Substructures
of}

\smallskip

\centerline{\tit Structures on ${\mtit \aleph_\omega}$}
\bigskip
\bigskip

$$\matrix{
\hfill \hbox{\bit Kecheng Liu}\hfill\cr
\hfill \hbox{Department of Mathematics}\hfill\cr
\hfill \hbox{University of California}\hfill\cr
\hfill \hbox{Irvine, CA 92717}\hfill\cr}
$$
\smallskip
\footnote{}{$^*$ Partially supported by the United
States-Israel Binational Science Foundation publ 484.}

\centerline{and}
\smallskip

$$
\matrix{
\hfill\hbox{\bit Saharon Shelah$^*$}\hfill\cr
\hfill\hbox{Institute of Mathematics}\hfill\cr
\hfill\hbox{The Hebrew University}\hfill\cr
\hfill\hbox{Jerusalem, Israel}\hfill\cr
\hfill\hbox{Department of Mathematics}\hfill\cr
\hfill\hbox{Rutgers University}\hfill\cr
\hfill\hbox{New Brunswick, NJ}\hfill\cr
%\hfill\hbox{Department of Mathematics}\hfill\cr
%\hfill\hbox{University of Wisconsin}\hfill\cr
%\hfill\hbox{Madison, WI}\hfill\cr
}
$$
\bigskip
Revisions made 3/29/96

\bigskip

\n {\bf 1. INTRODUCTION.}
\smallskip

Let $0<n^* < \omega$ and $f:X\to n^*+1$ be a function where
$X\subseteq \omega\backslash (n^*+1)$ is infinite.  Let's
consider the following set
$$
S_f=\{ x\subset\aleph_\omega :|x|\le\aleph_{n^*}\wedge
(\forall n\in X)cf(x\cap\alpha_n)=\aleph_{f(n)}\}.
$$

The question is whether $S_f$ is stationary in
$[\alpha_\omega]^{<\aleph_{n^*+1}}$.  The question was first
posed by Baumgartner in \cite{B7}.  By a standard result, the
above question can also be rephrased as certain transfer
property.  Namely, $S_f$ is stationary iff for any structure
$A=\la\aleph_\omega ,\ldots\ra$ there's a $B\prec A$ such
that $|B|=\aleph_{n^*}$ and for all $n\in X$ we have
$cf(B\cap\aleph_n)=\aleph_{f(n)}$.

Note that if $f$ is eventually constant then $S_f$ is always
stationary.  (see \cite{B7}).  Also, any reasonable ``finite
variation'' of $f$ will preserve the property, i.e., if
$n^*_1>n^*$ and $f_1:X\rightarrow n^*+1$ agrees with $f$ on
$\hbox{dom}(f)\backslash (n^*_1+1)$, then $S_{f_1}$ is
stationary provided that $S_f$ is.  So we are interested in
the case that $f$ is not eventually constant.  You may
wonder how strong the statement that $S_f$ is stationary is.
Magidor (unpublished, but close to \cite{Mg3}) has shown that if
$f:\omega\backslash 2\to 2$ assumes the values $0$ and $1$
alternatively and $S_f$ is stationary (here $n^*=1$), then
there's an inner model with infinitely many measurable
cardinals.

Liu has proven (in \cite{Liu}) that under the existence of huge
cardinals, it is consistent that $S_f$ is stationary
together with GCH for many $f$'s assuming (any) two
different values.  In this paper, we are going to prove a
few results concerning the above question, which are due to
the second author except for Proposition 3.1 and a small
remark of the first author improving Theorem 4.2.

Theorems here (3.1, 4.2) are
the first results for function $f$ with more than $2$
values gotten infinitely many times.  Also the present
results have relatively small consistent strength.  A version of
4.1 was proved in the mid eighties and forgotten and we thank
M. Gitik for reminding.  The main results are as follows:
\medskip

\n {\bf Theorem 3.1.}  {\it Assume $\sup
(pcf(\{\aleph_n:n<\omega\}))=\aleph_{\omega +n^*}$.  Let
$1<m^*<\omega$.  Let $I$ be the ideal of finite subsets of
$\omega$.  Let $\la A_i:1\le i\le n^*\ra$ be pairwise
disjoint subsets of $\omega\backslash (m^*+1)$ such that
$\prod_{k\in A_i}\aleph_k /I$ has true cofinality
$\aleph_{\omega +i}$ for $1\le i\le n^*$.  Let
$\la\kappa_i:1\le i\le n^*\ra$ be a sequence of uncountable
cardinals below $\aleph_{m^*+1}$.  Then the set
$$
S=\{ x\subseteq \alo :|x|\le \aleph_{m^*}\wedge (\forall
i) [ 1\le i\le n^*\to (\forall k\in A_i)cf(x\cap\aleph_\kappa
)=\kappa_i ] \}
$$
is stationary in $[\alo ]^{<\aleph_{m^*+1}}$.  }
\medskip

\n{\bf Theorem 4.1.}  {\it Assume $A\subseteq \omega$,
$0<n^* <\min (A)$ and for each $n\in A$ there is an
$\aleph_n$-complete filter $F_n$ on $\aleph_n$ such that
$F_n$ contains the cobounded subsets of $\aleph_n$ and the
second player has a winning strategy in the game
$GM_{\aleph_{n^*}}(F_n)$.
\smallskip
\n (See Definition 4.1.)
\medskip

Then the set
$$\eqalign{
S=&\left\{ x\subset\alo :|x|\le\aleph_{n^*}\;\hbox{and}\; (\forall
n\in A)[cf(x\cap\aleph_n)=\aleph_0]
\;\;\hbox{and}\right.\cr
&\left.(\forall n)[n<\omega \wedge
n^*<n\not\in A  \rightarrow
cf(x\cap\aleph_n)=\aleph_{n^*}]\right\}\cr}
$$
is stationary in $[\alo ]^{<\aleph_{n^*+1}}$.
\medskip

\n {\bf Theorem 4.2.} {\it Assume GCH.  Let $0\le m<
n^*<\omega$ and $E\subseteq\omega\backslash n^*$ be such
that for all $i\in E$, $i+1\not\in E$.  Let $\la n_i:i\in
E\ra$ with each $n_i <n^*$.  Suppose that for each $i\in E$,
there is an $\aleph_i$-complete filter $F_i$ on $\aleph_i$
containing all clubs of $\aleph_i$ such that $W_i=\{\alpha
<\aleph_i :cf(\alpha )=\aleph_{n_i}\}\in F^+_i$ and the
second player has a winning strategy in the game
$GM'_{\aleph_m}(F_i)$ (see Definition 4.1).  %  }
\medskip

Let $f:\omega\to n^*$ be the function defined by
$$
f(i)=\cases{n_i&if $\in E$,\cr
&\cr
m&if $n^* <i\not\in E$.\cr}
$$
Then the set $S'=\{ x\subseteq\alo :|x|\le \aleph_{n^*}$ and
$(\forall i> n^*)[cf (x\cap\aleph_i)=\aleph_{f(i)}]\}$ is
stationary in $[\alo ]^{<\aleph_{n^*+1}}$.     }  }
\medskip

In this paper, we concentrate on $\{\aleph_n :n\in\omega\}$,
but we can deal with other sets with natural changes (by the
referee request an explanation was added in 95 in the end of
the paper).  We implicitly assume that all models under
consideration have a countable language.
\bigskip

\n {\bf 2. PRELIMINARIES.}
\smallskip

Let's start with a standard result.  The proof will be
omitted.
\medskip

\n {\bf Proposition 2.1.} {\it Let $n^*<\omega$,
$X\subseteq\omega\backslash (n^*+1)$ and $f:X\to n^*+1$.
Consider the set $S=S_f$ as defined in \S1.  Let $\theta
>(\alo^{\aleph_{n^*}})^+$ be a regular cardinal.  Let
$M\prec \la \HH_\theta, \in , S, \triangleleft , \ldots \ra$,
$M\supseteq\alo +1$ and $|M|=\alo$, where $\triangleleft$ is
a well-ordering of $\HH_\theta$.  Then the following are
equivalent:
\smallskip
\item{1.}  $S$ is stationary in $[\alo ]^{<\aleph_{n^*+1}}$.
\smallskip

\item{2.} For any structure $A=\la\alo ,\ldots \ra$ with a
countable language, there is a $B\prec A$ such that
$|B|=\aleph_{n^*}$ and $B\in S$.
\smallskip

\item{3.}  There is $N\prec M$ such that $|N|=\aleph_{n^*}$
and $\forall n\in X$, $cf(N\cap \aleph_n)=\aleph_{f(n)}$.}
\medskip

\n  {\bf Lemma 2.1.} {\it Let $\kappa <\mu <\lambda$ be
regular cardinals.  Let $A=\la \HH_\lambda, \in,
\triangleleft, \kappa, \mu ,\ldots\ra$ be a structure of a
countable language on $\HH_\lambda$ with skolem functions
closed under composition and $\triangleleft$ a well-ordering
of $\HH_\lambda$.  If $B\prec A$ and $X\subseteq\kappa$ and
$B'=sk^A(B\cup X)$, then $\sup (B'\cap\mu )=\sup (B\cap \mu
)$.}
\medskip

\n {\bf Proof.}  It's clear that the lemma holds if $\sup
(B\cap\mu )=\mu$.  So we assume $\sup (B\cap\mu )<\mu$.

It's clear that $\sup (B\cap\mu )\le\sup (B'\cap\mu )$.  Now
suppose $\alpha\in B'\cap\mu$.  WLOG, assume $\alpha =\tau
(b,x_0)$ for some $b\in B$, $x_0\in X$ and some skolem
function $\tau$.  Define $f:\kappa\to\mu$ by letting
$f(x)=\tau (b,x)$ if $\tau (b,x)<\mu$ and $f(x)=0$
otherwise.  Then $f$ is definable from $b$ in $B$.  So $f\in
B$.  Let $\delta =\sup (f''\kappa )$.  Then $\delta\in B\cap
\mu$ and $B'\models \delta =\sup (f''\kappa )$.  So
$B'\models \alpha =f(x_0)\le\delta$.  Hence $\alpha\le\delta
<\sup (B\cap\mu )$.  Therefore, $\sup (B\cap\mu )\ge\sup
(B'\cap\mu )$.  This completes the proof of the
lemma.~$\square$
\medskip

\n {\bf Lemma 2.2.} {\it Let $n\in\omega$ and $X_0$,
$S_1,\ldots , X_{n-1}$ be disjoint subsets of $\omega$.  Let
$X=\bigcup_{i<n}X_i$ and $f$ be a function from $X$ to
$\omega$ such that $f$ is constant on $X_i$ for each $i<n$
and the constant values of $f$ on different $X_i$'s are
distinct.  Let $A=\la \alo ,\ldots\ra$ be an algebra on
$\alo$.  Let $B\prec A$ be such that $(\forall i<n)$
$(\forall m\in X_i)$ $[cf(B\cap\aleph_m)=\aleph_{f(m)}]$.
Let $k=\max (f''X)$, $A_{i_0}=f^{-1}\{ k\}$ and $\ell <k$.
Then for any $n^*$ such that $n^*\ge \ell$ and $n^*>\max
(f''X\backslash\{ k\})$, and $j<\omega$ such that $|B\cap
\aleph_j|=\aleph_{n^*}$, there is $B'\prec B$ such that
\smallskip
\item{(1)} $|B'|=\aleph_{n^*}$ and
$B'\cap\aleph_m=B\cap\aleph_m$ for $m \leq j$;
\smallskip
\item{(2)} $(\forall m\in X_{i_0})$ $[m>j\to
cf(B'\cap\aleph_m)=\aleph_\ell]$;
\smallskip
\item{(3)} $(\forall m\in X\backslash X_{i_0})$ $[m>j\to cf
(B'\cap\aleph_m)=cf(B\cap \aleph_m)]$.}
\medskip

\n {\bf Proof.}  For each $i\not= i_0$, for each $m\in X_i$, let
$a_m$ be a cofinal subset of $B\cap\aleph_m$ with order type
$\aleph_{f(i)}$.  Now we can build a sequence $\la B_\alpha
:\alpha <\aleph_\ell\ra$ such that
\smallskip
\item{(1)} $\cup\{ a_m:m\in X\backslash X_{i_0}\}\cup
(B\cap\aleph_j)\subseteq B_0$;
\smallskip
\item{(2)} $\forall m\in X_{i_0} [ \sup
(B_\alpha\cap\aleph_m)<\sup (B_{\alpha +1}\cap\aleph_m)$];
\smallskip
\item{(3)} $B_\alpha \prec B_{\alpha +1}\prec B$ and
$|B_\alpha |=\aleph_{n^*}$.
\smallskip

\n The construction is obvious.  Now let $B'=\bigcup_{\alpha
<\aleph_\ell}B_\alpha$.  It is clear that $B'$ is as
required.~$\square$
\medskip

\n {\bf 3. AN APPLICATION OF PCF THEORY.}
\smallskip

We are going to prove the following theorem using pcf theory
(see \cite{Sh:g}).
\medskip

\n {\bf Theorem 3.1.} {\it Assume $\max (pcf
(\{\aleph_n:n<\omega\}))=\aleph_{\omega +n^*}$.  Let
$1<m^*<\omega$.  Let $I$ be the ideal of finite subsets of
$\omega$.  Let $\la A_i:1\le i\le n^*\ra$ be a sequence of
pairwise disjoint subsets of $\omega\backslash (m^*+1)$ such
that $\prod_{k\in A_i}\aleph_k/I$ has true cofinality
$\aleph_{\omega +i}$ for $1\le i\le n^*$.  Let
$\la\kappa_i:1\le i\le n^* \ra$ be a sequence of uncountable
cardinals below $\aleph_{m^*+1}$.  Then the set
$$
S=\{ x\subseteq\aleph_\omega
:|x|\le\aleph_{m^*}\wedge(\forall i)(1\le i\le n^*  \to
(\forall k\in A_i)[cf(x\cap\aleph_k)=\kappa_i])\}
$$
is stationary in $[\alo ]^{<\aleph_{m^*+1}}$.}

\bigskip

\n {\bf Remarks.}
\smallskip
\item{1.} Using Lemma 2.2, Lemma 2.1 and this theorem, we
can show that for any given $1\le n\in\omega$ and $\la
\kappa_i:1\le i\le n^*\ra$ with cardinals
$\kappa_i\le\aleph_n$ for $1\le i\le n^*$, we have that the
set
$$
S'=\{ x\subset\alo :|x|\le\aleph_n\;\hbox{and}\; (\forall
i)(1\le i\le n^*  \to \forall k\in
A_i[cf(x\cap\aleph_i)=\kappa_i])\}
$$
is stationary in $[\alo ]^{<\aleph_{n+1}}$.
\smallskip
\item{2.} Notice that in order for the theorem not to be
trivial, we assume $n^*>1$ and therefore GCH fails at
$\alo$.
\smallskip
\item{3.} If $pp\alo \mo \sup (pcf
(\{\aleph_m:m<\omega\}))>\aleph_{\omega +n^*}$, no harm is
done since we can use Levy collapse to collapse $pp\alo$ to
$\aleph_{\omega +n^*}$ and no new subset of $\alo$ is added.
\smallskip
\item{4.} The theorem can be generalized to other singular
cardinals.  Also, we can use other regular cardinals in
$(\alo ,pp\alo )$ in the proof of the theorem.
\smallskip
\item{5.} Consistency results giving the assumptions are
well know, starting with \cite{Mg}; see history and references
on this in \cite{Sh:g}.
\smallskip

\n The proof of Theorem 3.1 uses the following lemma:
\medskip

\n {\bf Lemma 3.1.}  {\it For each $1\le i\le n^*$, there is
a sequence $\vec C^i \mo \la C^i_\alpha :\alpha <\aleph_{\omega +1}\ra$ such
that
\smallskip
\item{(1)} $\forall \alpha C^i_\alpha\subseteq\alpha$ and
o.t.$(C^i_\alpha )\le\kappa_i$
\smallskip
\item{(2)} $\beta\in C^i_\alpha$ implies $C^i_\beta
=C^i_\alpha\cap\beta$
\smallskip
\item{(3)} $S_i \mo \{\alpha <\aleph_{\omega +i}:cf(\alpha )=\kappa_i\;\;
\hbox{and}\;\; \alpha=\sup (C^i_\alpha )\}$ is stationary in
$\aleph_{\omega +i}$.}
\medskip

\n {\bf Remarks.}  Note that the $C^i_\alpha$'s are not
necessarily closed.
\medskip

For a proof of Lemma 3.1, see \cite[4.4]{Sh:351} for successor
of regular cardinals and in general \cite[1.5]{Sh:420} which
rely on \cite[4.4]{Sh:351}.

For each $1\le i\le n^*$, let $S_i$ and $\vec C^i=\la
C^i_\alpha :\alpha <\aleph_{\omega +i}\ra$ be as in Lemma
3.1.  We now proceed to prove Theorem 3.1.
\medskip

\n {\bf Proof of Theorem 3.1.}  For $1\le i\le n^*$, let
$\vec f^i=\la f^i_\alpha :\alpha <\aleph_{\omega
+i}\ra\subseteq\Pi_{k\in A_i}\aleph_k$ be increasing and
cofinal in $\Pi_{k\in A_i}\aleph_k/I$.  Let
$\lambda>>pp\alo$ be a regular cardinal.  Let's consider the
structure $A=\la\HH_\lambda , \in ,\triangleleft,\ldots\ra$
with skolem functions closed under composition, where
$\triangleleft$ is a well-ordering on $\HH_\lambda$.  We
define $X_{\vec\alpha}, N_{\vec\alpha}$ by induction on
$\vec \alpha_{n^*}$ as follows.  Let $x=\{\vec f^i, A_i, \vec
C^i\}_{1\le i\le n^*}$.  For each $\vec\alpha
=\la\alpha_i:1\le i\le n^*\ra\in\Pi_{1\le i\le
n^*}\aleph_{\omega +i}$, let
%%\bye
$$\eqalign{
X_{\vec\alpha}= & \{\gamma
:\gamma <\aleph_{n^*}\} \cup (\bigcup_{1\le i\le
n^*}C^i_{\alpha_i})\cup\{
N_{\vec\beta}:\vec\beta\in\prod^n_{i=1}C^i_{\alpha_i}\}\cr
& ( \hbox {hence}\;\;    C^i_\xi \subseteq X_{\vec\alpha}
\;\; \hbox {if} :   \;\; 1\le i\le n^*\wedge\xi \in
C^i_{\alpha_i} )  \cr}
$$
and let $N_{\vec\alpha}=sk^A(X_{\vec\alpha})$.  Note that
$|N_{\vec\alpha}|=\aleph_n$ and
$\vec\alpha\in\prod_{i=1}^{n^*}C^i_{\beta_i}\Rightarrow
N_{\vec\alpha}\prec N_{\vec\beta}$ and $N_{\vec\alpha}\in
N_{\vec\beta}$.
\medskip

\n {\bf Claim.}  There is $\vec\delta =\la\delta_i:1\le i\le
n^*\ra\in\Pi_{1\le i\le n^*}S_i$ such that all $1\le i\le
n^*$
\smallskip
\item{1)} For all $\vec\alpha\in\Pi_{1\le i\le
n^*}C^i_{\delta_i}$, we have $\sup (N_{\vec\alpha}\cap\aleph_{\omega
+i})<\delta_i$ for all $1\le i <  n^*$.
\smallskip

\item{2)} $\sup (N_{\vec\delta}\cap\aleph_{\omega
+i})=\delta_i$.
\smallskip

\item{3)} For some $n_i<\omega$, $\{ f^i_\alpha
(k):\alpha\in C^i_{\delta_i}\}$ is cofinal in
$N_{\vec\delta}\cap\aleph_k$ for all $k\in A_i\backslash
n_i$.
\smallskip

\item{4)} For some $m_i\ge n_i$,
$cf(N_{\vec\delta}\cap\aleph_k)=cf(\delta_i)=\kappa_i$ for
all $k\in A_i\backslash m_i$.
\medskip

\n {\bf Proof of Claim.}  We first construct $\vec\delta$ as
required by part 1) of the claim.  The construction is as
follows:  Let $E_{n^*}=\{\delta <\aleph_{\omega
+n^*}:\forall\vec\alpha\in\Pi_{1\le i\le n^*}\aleph_{\omega
+i}$, if $\alpha_{n^*}<\delta$, then $\sup
(N_{\vec\alpha}\cap\aleph_{\omega +n^*})<\delta\}$.  Then
$E_{n^*}$ is clearly closed unbounded in $\aleph_{\omega
+n^*}$.  Since $S_{n^*}$ is stationary in $\aleph_{\omega
+n^*}$ we have $S_{n^*}\cap E_{n^*}\not=\emptyset$.  Pick
some $\delta_{n^*}\in S_{n^*}\cap E_{n^*}$.

Suppose we have defined $\delta_j$ for $n^*\ge j>i$.  We now
define $\delta_i$.  Let $E_i=\{\delta <\aleph_{\omega
+i}:\forall\vec\alpha\in\Pi_{1\le \ell\le n^*}\aleph_{\omega
+\ell}$ if $\alpha_i<\delta$ and $\alpha_j=\delta_j$ for all
$i<j\le n^*$ then $\sup (N_{\vec\alpha}\cap\aleph_{\omega
+i})<\delta\}$.  It's easy to see that $E_i$ is closed
unbounded in $\aleph_{\omega + i}$.  So we can find
$\delta_i\in S_i\cap E_i$ since $S_i$ is stationary in
$\aleph_{\omega +i}$.

We now show $\vec\delta$ satisfies clause (1).  Let
$\vec\alpha\in\Pi_{1\le i\le n^*}C^i_{\delta_i}$.  Let $1\le
i\le n^*$ be fixed, and we want to show that $\sup
(N_{\vec\alpha}\cap\aleph_{\omega +i})<\delta_i$.  Consider
$\vec\beta =\la\alpha_1,\ldots \alpha_i,
\delta_{i+1},\ldots\delta_{n^*}\ra$.  By the choice of
$\delta_i$, we have that $\sup (N_{\vec
\beta}\cap\aleph_{\omega +i})<\delta_i$.  Since
$\vec\alpha \in \prod_{1\leq
k \leq n^*} C^k_{\delta_k}$ clearly
$X_{\vec \alpha} \subseteq X_{\vec
\beta}$.  So $\sup (N_{\vec\alpha} \cap \aleph_{\omega +i}) \leq \sup
(N_{\vec \beta} \cap \aleph_{\omega +i} ) < \delta_i$.

Let's now prove clause (2) of the claim.  Fix $1\le i\le
n^*$.  Let $\beta\in N_{\vec\delta}\cap\aleph_{\omega +i}$.
Then $\beta=\tau (\vec y)$ for some $\vec y\in
[X_{\vec\delta}]^{<\omega}$ and some skolem function $\tau$.
We need to show $\beta <\delta_i$.  Since $\delta_j=\sup
(C^j_{\delta_j})$ and $C^j_{\delta_j}$ has no last element
for $1\le j\le n^*$, there is $\vec\alpha\in\Pi_{1\le j\le
n^*}C^j_{\delta_j}$ such that $\vec y\in
[X_{\vec\alpha}]^{<\omega}$.  But then $\beta\in
N_{\vec\alpha}\cap\aleph_{\omega +i}$.  By clause (1) we
have $\beta <\delta_i$ so $\sup (N_{\vec\alpha} \cap \aleph_{\omega + 1})
\leq \delta_i$ and by the previous paragraph we get equality.

We show clause (3) by contradiction.  Assume that (3) fails.
So there is an unbounded set $b\subseteq A_i$ such that
$\forall k\in b\exists\beta_k\in
N_{\vec\delta}\cap\aleph_k[ \sup (\{ f^i_\xi (k):\xi\in
C^i_{\delta_i}\})<\beta_k]$.  Fix such $\beta_k$ for each
$k\in b$.  Since for each $1\le j\le n^*$ the set
$C^j_{\delta_j}$ has order whose cofinality is uncountable,
there is $\vec\alpha\in\Pi_{1\le\ell \le
n^*}C^\ell_{\delta_\ell}$ such that $( \forall k\in
b) \beta_k\in N_{\vec\alpha}$.  Clearly $k\in
B\Rightarrow\beta_k < \sup (N_{\vec\alpha}\cap\aleph_k)$ hence
$\la\beta_k:k\in b\ra <\la\sup
(N_{\vec\alpha}\cap\aleph_k):k\in b\ra$.  The later belong
to $\Pi_{k\in b}\aleph_k$ hence for some $\xi < \aleph_{\omega
+i}$ we have $\la\sup (N_{\vec\alpha}\cap\aleph_k):K\in
b\ra\le^*f^i_\xi$, where $f\le^* g$ means $f(k)\le g(k)$ for
all but finitely many $k$'s.  Since $N_{\vec\alpha}\in
N_{\vec\delta}$, clearly $\la\sup
(N_{\vec\alpha}\cap\aleph_k ) :k\in b)$ belongs to
$N_{\vec\delta}$ and also $\la f^i_\xi :\xi
<\aleph_{\omega +i}\ra$ belongs to $X_{\vec\delta}$ hence to
$N_{\vec\delta}$, so wlog $\xi\in
N_{\vec\delta}\cap\aleph_{\omega +i}$.  Now we can replace
$\xi$ by any $\xi'\in (\xi ,\aleph_{\omega +1})$ and
$C^i_{\delta_i}$ is unbounded in
$N_{\vec\delta}\cap\aleph_{\omega +i}$, so wlog $\xi\in
C^i_{\delta_i}$ hence $k\in b\Rightarrow f_\xi (k)\in
N_{\vec\delta}\cap\aleph_{\omega +i}$ hence is $<
\beta_k$.  This is clearly absurd.

Finally, let's prove clause (4) again by contradiction.
Suppose clause (4) is not true.  Then there is an unbounded
set $b\subseteq A_i\backslash n_i$ such that $cf(\sup\{
f^i_\xi (k):\xi\in C^i_{\delta_i}\} )<\kappa_i$ for some
$k\in b$.  So as o.t.$(C^i_{\delta_i})=\kappa_i$ for each
$k\in b$, there is $\xi_k\in C^i_{\delta_i}$ such that for
all $\xi\in C^i_{\delta_i}$ with $\xi\ge\xi_k$ we have:
$\sup\{ f^i_\zeta (k):\zeta\in C^i_{\delta_i}\}=\sup\{ f^i_\zeta
(k):\zeta\in C^i_{\delta_i}\cap\xi\}$.  Let $\beta\in
C^i_{\delta_i}$ be such that $\beta >\sup (\{\xi_k :k\in
b\})$.  For each $k\in A_i$, let $\mu_k=\sup (\{ f^i_\zeta
(k):\zeta \in C^i_\beta\} )$.  Then $\la\mu_k:k\in A_i\ra\in
N_{\vec\delta}\cap\Pi_{k\in A_i}\aleph_k$ since
$C^i_\beta\in N_{\vec\delta}$.  So as above there is
$\beta'\in N_{\vec\delta}\cap\aleph_{\omega +i}$ such that
$\la\mu_k:k\in A_i\ra\le^* f^i_{\beta'}$.  So we have
$\la\sup (N_{\vec\delta}\cap\aleph_k ):k\in b\ra
=\la\mu_k:k\in B\ra\le^*f^i_{\beta'} \upharpoonright b$ which
contradicts to $f^i_{\beta'}\in N_{\vec\delta}\cap\Pi_{k\in
A_i}\aleph_k$ (the initiated could have used ``wlog $\vec
f^i$ obeys $\vec C^i$'').  This completes the proof of the
claim.

Now, Theorem 3.1 follows from the claim, Proposition 2.1 and
Lemma 2.1.~$\square$
\medskip

By the remarks following Theorem 3.1, the theorem is not
trivial only when $pp\alo >\aleph_{\omega +1}$.  In
particular GCH does not hold at $\alo$.  But using the
following observation, we can make GCH hold at $\alo$ by
collapsing $2^{\alo}$ and still have the desired conclusion
in the forcing extension, i.e., the set $S$ in Theorem 3.1
is still stationary in the generic extension.  (So by
well-known consistency results we can even have GCH.)
\medskip

\n {\bf Proposition 3.1.}  {\it Let $P$ be an
$<\aleph_{\omega +1}$-closed forcing notion.  Suppose $S$ is
stationary in $[\alo ]^{<\aleph_{n+1}}$ in $V$.  Then
$V^P \Vdash$ ``$S$ is stationary in $[\alo
]^{<\aleph_{n+1}}$''.}
\medskip

\n {\bf Proof.}  It suffices to show that in $V^p$, for any
given structure $A=\la\alo ,\ldots\ra$ of a countable
language, there is $B\prec A$ such that $|B|=\aleph_n$ and
$B\in S$.

Let $p\in P$ for that $\dot A=\la\alo ,(\dot
f_i)_{i\in\omega}\ra$ is a structure on $\alo$ with skolem
functions $\dot f_i$ closed under compositions.

Since $P$ is $( <\aleph_{\omega +1} ) $-closed, we can find $\la
p_\alpha :\alpha <\alo\ra$ (such that $p_\alpha$ is stronger
than $p$ and $p_\beta$ for $\beta <\alpha$ and $f'_i$ in $V$
such that for each $i$, for any $\vec\alpha\in [\alo
]^{<\omega}$ there is a $\beta$ such that $p_\beta\Vdash\dot
f_i(\vec\alpha )=f'_i(\vec\alpha )$ whenever $\vec\alpha\in
\hbox{dom}(\dot f)$.

Consider $A'=\la\alo ,(f'_i)_{i\in\omega}\ra$ in $V$.  Let
$p$ be such that $p$ is stronger than $p_\alpha$ for all
$\alpha <\alo$.  Let $B\prec A'$ be such that $|B|=\aleph_n$
and $B\in S$.  Sut then $p\Vdash B\prec\dot A$ since
$p\Vdash f'_i=\dot f_i$.  This is as required.~$\square$
\medskip

\n {\bf 4. APPLICATIONS OF LARGE IDEALS.}
\smallskip

In this section, we prove two results under the existence of
large ideals (on the $\aleph_n$'s).  Before we state our
results, we need some terminology.

%%%%%%%%%%%

\bigskip

\noindent {\bf Definition 4.1.}  {\it Let $\kappa > \lambda$ be cardinals.
Let $D$ be a filter on $\kappa$.

\noindent (1)  We define the game $\underline {GM_\lambda (D)}$ as follows:
the game lasts $\lambda$ moves.  At $\xi^{th}$ move, the first player
chooses a subset $A_\xi$ of $\kappa$ such that  $A_\xi \subseteq \cap_{\eta <
\xi} B_\eta$, and if $\cap_{\eta < \xi} B_\eta \not = \emptyset \mod (D)$
then $A_\xi \not = \emptyset \mod (D)$.  The second player chooses a sbuset
$B_\xi$ of $A_\xi$ with $B_\xi \not = \empty \mod (D)$.

A player without a legal move loses the game immediately.  (Note this can
only happen to the second player.)  If the game lasts for $\lambda$ moves,
the second player wins if $\cap_{\xi < \lambda} B_\xi$ is unbounded.

\noindent (2)  Let's also assume that $D$ is $\kappa$-complete.  We now
define the ``cut-and-choose'' game $\underline {GM'_{\lambda,\kappa} (D) }$
of length $\lambda$:  at the $0^{th}$ move, the first player chooses a set
$A_0 \not = \emptyset \mod (D)$ and then partitions $A_0$ into less than
$\kappa$ parts; the second player chooses one of the parts, say $B_0 \in
D^+$.  At the $\xi^{th}$ move for $\xi > 0$, the first player partitions the
set $\cap_{\eta < \xi} B_\eta$ into less than $\kappa$ parts, and the second
player chooses one of the parts, call it $B_\xi$ such that $B_\xi \in D^+$.

The winning conditions for each player is exactly as in the game defined in
part (1) above.

}

Let's first prove the following theorem:

\bigskip

\noindent {\bf Theorem 4.1.}  {\it Assume $A \subseteq \omega, 0 < n^* < \min
(A)$ and for each $n \in A$ there is an $\aleph_n$ such that $F_n$ contains
the cobounded subsets of $\alpha_n$ and the second player has a winning
strategy in the game $GM_{\alpha_{n^*}} (F_n)$.

Then the set
$$
S = \{ x \subset \aleph_\omega : |x| \leq \alpha_{n^*} ,
  (\forall n \in A) cf (x \cap \aleph_n) =
     \aleph_0 \;\;\hbox{and} \;\;
         (\forall n) [ \forall n^* < n \not \in A \rightarrow
         cf (x \cap \aleph_n) =
         \aleph_{n^*} ] \}
$$
is stationary in $[\aleph_\omega ]^{{<\aleph}_{n^* + 1}} $.

}

\bigskip

\noindent {\bf Remarks.}

\noindent (1)  Before we prove Theorem 4.1, we would like to see how we can
get the filters as required in the hypothesis of the theorem.  Magidor (see
\cite{Mg1}) has shown the consistency of the existence of the filters.  Also,
Laver has proved that if we collapse a measurable cardinal $\kappa$ to
some $\aleph_n$, then in the generic extension, there is a normal ideal on
$F_n$ on $\aleph_n$ such that $\PP (\aleph_n) / F_n$ has a $<
\aleph_{n-1}$-closed dense subset.  Therefore, if there are infinitely many
measurable cardinals, say $\la \kappa_n : n < \omega \ra$,
$A = \{   m_m | n < \omega \}, n^* < m_n < n, m_n + 1 < m_{n + 1} $
we can Levy collapse
each $\kappa_n$ to $\aleph_{n^* + m_n}$ to get the normal filters as required
in the hypothesis of Theorem 4.1.

\noindent {(2)}  We can also use (in the assumption of Theorem 4.1) the
games $GM'_{\aleph_{n^*} , \aleph_n }   (F_n)$ in place of $GM_{\aleph_{n^*}}
(F_n)$.  We can weaken it further using for $n \in A$, the following game
for $F_n$ (see \cite{Sh:250}) (see better \cite[Ch XIV]{Sh:f}: in the
$\xi^{th}$ move, player one choose $m_\xi \in \omega \backslash A, n^*
< m_\xi < n$ and $F_\xi : \aleph_n \rightarrow \aleph_{m_\xi}$ and
player two has to choose $B_\xi \subseteq \cap_{\zeta < \xi} B_\xi
\cap A_0 $ such that the range of $f_\xi 
\upharpoonright B_\xi$ is bounded in $\aleph_{m_\xi}$ and $B_\xi \not =
\emptyset \mod F_n$; in the $0^{th}$ move player one also choose $A_0
\subseteq \aleph_{m_\xi}, A_0 \not =  \emptyset \mod F_n$.

If $F$ is a filter on
$\lambda$, and the cardinal player two choose is from $S$,
and a play last $\theta$ moves we call the game $GM_\theta'
(F,S)$.

In order to prove Theorem 4.1, let's consider tagged trees of the form $\la
T,I \ra$, which by definition means that

\smallskip

\item {1.}  $T \subset [ON]^{<\omega}$ is a tree, i.e. $T$ consists of
finite sequences of ordinals closed under initial segments.

\item {2.}  $\II = \la I_\sigma : \sigma \in T\ra$ is such that for each
$\sigma \in T, I_\sigma$ is an ideal on $Suc_T(\sigma)$ which is the set of
immediate successors of $\sigma$ in $T$.  Also, $I_\sigma$ can be thought of
as an ideal on $\{ \alpha : \sigma^{\la \alpha \ra} \in T\}$.

If $T_1$ is a subtree of $T$, we can view $\la T_1, \II'\ra$ as a tagged tree
with $\II' = \la I_\sigma |Suc_{T_1} (\sigma ) : \sigma \in T_2 \ra$.  By
abuse of notation, we still denote it by $\la T_1, \II\ra$.
If the family $\II$
of ideals is clear from context, we will simply say $T$ is a tagged tree
without mentioning $\II$ explicitly.

For $X \subseteq T$, let $T[X] = \{ \sigma \in T: \exists \eta \in X(\eta
\leq_T \sigma \vee \sigma \leq_T \eta)\}$.  Clearly $T[X]$ is a subtree of
$T$.  The following lemma is from \cite{RuSh:117} or \cite[Ch X]{Sh:b}
or \cite{Sh:f} and we will not give the proof here.

\bigskip

\noindent {\bf Lemma 4.1.}  {\it Let $\la T,\II  \ra$ be a tagged tree
such that
for each $\sigma \in T, I_\sigma$ is a proper ideal such that $Suc_T (\sigma) \not
\in I_\sigma$.  Let $\lambda$ be a
regular, uncountable cardinal and for every $\sigma \in T, I_\sigma$ is
$\lambda$-indecomposable, i.e. if $A \subseteq Suc_T (\sigma ) A\not =
\emptyset \mod I_\sigma$ and $f:A \rightarrow \lambda$ then for some $\zeta
< \lambda$ we have $ \{ x \in A: f(x) < \zeta\} \not =
\emptyset \mod I_\sigma$.  (This
holds if for every $\sigma \in T$, $I_\sigma$ is a $\lambda^+$ complete
ideal or  $ | \Suc_T (\sigma ) | < \lambda)$.

\item {(*)}  for every function $F : T \rightarrow \lambda$, there is a
subtree $T_1$ of $T$ such that for all $\sigma \in T_1$, $Suc_{T_1} (\sigma)
\not\in I_\sigma$ and $\Sup (F''T_1) < \lambda$.

}

\bigskip

We now proceed to prove Theorem 4.1.

\bigskip

{\bf Proof of Theorem 4.1.}  Let $\la m_i : i < \omega \ra$ be such that each
$m_i \in A$ and for each $m \in A$ there are infinitely many $i$ with $m =
m_i$.  Let $\TT =
\left \{  \right. T : T$ is a subtree of
$\bigcup_{l<\omega} \prod_{i<l} \aleph_{m_i}$ and $\forall\eta \in T$ with
$lh(\eta) = i$ we have ,
$\left \{   \right. \alpha < \aleph_{m_i} : \eta^{\hat{} } \la \alpha
\ra \in T
\left. \right \}
\in F^+_{m_i}
\left. \right \}$, where $lh(\eta )$ means the length of the finite sequence
$\eta$.

Note that each $T \in \TT$ can be considered as a tagged tree where for each
$\sigma \in T$, the associated ideal $I_\sigma$ is just the dual ideal to
the filter $F_{m_{lh(\sigma)}}$.

Suppose $A = \la \aleph_{\omega , \ldots} \ra$ is an arbitrary structure on
$\aleph_\omega$.  We are going to find a $B \prec A$ such that $|B | =
\aleph_{n^*}$ and for each $n^* < m < \omega$ if $m \in A$ then $cf(B \cap
\aleph_m) = \aleph_0;$ if $m \not \in A$ then $cf(B \cap \aleph_m) =
\aleph_{n^*}$.  This is enough to prove the theorem by Proposition 2.1.

By induction on $\xi < \aleph_{n^*}$, we are going to build
$T_\xi, \la \alpha_{\xi ,m} : n^* < m \not\in A\ra$,
$\la A_{\xi ,\eta}, B_{\xi ,\eta} : \eta \in T_{\xi + 1} \ra$ and
$\la N_{\xi , \eta} : \eta \in T_\xi\ra$ such that

\item {1.}  $T_\xi \in \TT$ and for any $\xi < \xi', T_{\xi'} \subseteq
T_\xi$.

\item {2.}  For $\eta \in T_{\xi + 1}, B_{\xi ,\eta} =
   \{ \alpha < \aleph_{m_{lh(\eta )} } : \eta^{\hat{}} < \alpha >\in T_{\xi +
   1} \}$.  Furthermore, $\la A_{\xi',\eta} , B_{\xi ',\eta} : \xi ' \leq
   \xi \ra$ is an initial segment of a play of the game $GM_{\aleph_{n^*}}
   (F_{m_{lh(\eta)} } )$ with the second player following his winning
   strategy.

\item {3.}  If $\xi$ is a limit ordinal, $T_\xi = \cap_{\xi ' < \xi} T_{\xi
'}$.

\item {4.}  For $\eta \in T_\xi, N_{\xi ,\eta} = sk^A(ran(\eta ) \cap \{
\alpha_{\xi ',m} :
    \xi ' < \xi \;\; \hbox {and} \;\; n^* < m \not \in A\} ) .$

\item {5.}  For $\eta \in T_{\xi + 1}$ and $n^* < m \not \in A,$
we have $\Sup (N_{\xi
+ 1,\eta} \cap \aleph_m ) < \alpha_{\xi + 1,m}.$

\item {6.}  For each $n^* < m \not \in A, \la \alpha_{\xi ,m} : \xi <
\aleph_{n^*} \ra$ is an increasing sequence of ordinals in $\aleph_m.$

Take any $T_0 \in \TT$ to start with.  For $\xi$ limit, let $T_\xi =
\cap_{\xi ' < \xi} T_{\xi '}.$  Since the second player has a winning
strategy in the game $GM_{\aleph_{n^*} } (F_m)$ for each $m \in A, T_\xi \in
\TT$ for $\xi$ limit.

If $\xi = 0$, we let $\alpha_{\xi , m} = 0$ for $n^* < m \not \in A$.  If
$\xi$ is limit, we let $\alpha_{\xi , m} = \Sup(\{ \alpha_{\xi ',m} : \xi ' <
\xi \})$ for $n^* < m \not \in A$.

Suppose $T_\xi$ and $\la \alpha_{\xi , m} : n^* < m \not \in A\ra$ have been
constructed.  We now construct $T_{\xi + 1},
\la \alpha_{\xi + 1,m} : n^* < m \not \in A\ra$ and
$\la A_{\xi ,n}, B_{\xi , \eta} : \eta \in T_{\xi + 1}\ra$.

Let $\la k_i : i < \omega \ra$ be an enumeration of $\{ m \not \in A: n^* <
m < \omega \}$.  We will define $\la T'_i : i < \omega\ra,
\la A_{i,\eta}, B_{i,\eta} : i < \omega \ra$ for $\eta \in T'_{i+1}$ and
$\la \alpha_i : i < \omega \ra$ by induction on $i$ such that

\item {1.}  for each $i, T'_{i+1} \subseteq T'_i \in \TT$ and $\Sup (N_{\xi +
1,\eta} \cap \aleph_{k_i}) < \alpha_i$ for $\eta \in T'_{i+1}$;

\item {2.}  $\la A_{i,\eta}, B_{i,\eta} : i \in \omega \ra$ is an initial
segment of the play of the game $GM_{\aleph_{n^*}} (F_{lh(\eta )} )$ with
player two following his winning strategy.

Let $T'_0 = T_\xi$.  Suppose we have defined $T'_i, \alpha_{i-1}$ and
$A_{i-1,\eta}, B_{i-1,\eta}$ for $\eta \in T'_i.$  Consider the function
$F: T'_i \rightarrow \aleph_{k_i}$ defined by
$F(\eta ) = \Sup (N_{\xi + 1,\eta} \cap \aleph_{k_i})$.  Then $F$ has a value
$< \aleph_{k_i}$.  Since $k_i \not \in A$, we have $m_i \not = k_i$, so
$F_{m_i}$ is $\aleph_{m_i}$-complete on a set of cardinality $\aleph_{m_i}$
so the assumptions of Lemma 4.1 holds.  Hence there is $T''_{i+1} \subseteq
T'_i$ such that $T''_{i+1} \in \TT$ and $\Sup (F'' T'_{i+1}) < \aleph_{k_i}$
by Lemma 4.1.  Let $A_{i,\eta} = \Suc_{T''_{i+1}} (\eta)$ for $\eta \in
T''_{i+1}.$  Let $B_{i,\eta}$ be the move of the second player following his
winning strategy in the game $GM_{\aleph_{n^*}} (F_{lh(\eta)} )$.  Let
$T'_{i+1} \in \TT$ be such that $B_{i,\eta} = Suc_T(\eta)$ for each $\eta
\in T$.  Let $\alpha_i$ be such that $\Sup (F'' T'_{i+1} ) < \alpha_i <
\aleph_{k_i}.$

Now, let $T'_{\xi + 1} = \bigcap_{i<\omega} T'_i$ and $\alpha_{\xi + 1,k_i}
= \alpha_i$.  Since $\la A_{i,\eta}, B_{i,\eta} : i \in \omega\ra$ is an
initial segment of the play of the game $GM_{\aleph_{n^*} } (F_{lh(\eta)} )$
with player two following his winning strategy, we have that $T'_{\xi + 1}
\in \TT$.  For each $\eta \in T'_{\xi + 1}$, let $A_{\xi ,\eta} =
Suc_{R'_{\xi + 1} } (\eta)$.  Note this is a legal move for the first player.
Now player two chooses $B_{\xi ,\eta}$ according to his winning strategy.
Let $T_{\xi + 1} \in \TT$ be such that $B_{i,\eta} = Suc_{T_{\xi + 1}}
(\eta)$ for each $\eta \in T_{\xi + 1}$.  This completes the construction as
required.  (Alternatively demand that in $\la k_i: i < \omega\ra$, each
$m \in \{ m < \omega: m \not\in A, n^* < m < \omega \}$ appear
$\omega$ many times and if
$\xi = i \mod \omega$, take care only of $\aleph_{k_i}$.)

Finally, let $T_{\aleph_{n^*} } = \bigcap_{\xi < \aleph_{n^*}}$.  Since
$\la A_{\xi ,\eta}, B_{\xi , \eta} : \xi < \aleph_{n^*} \ra$ is a play of
the game $GM_{\aleph_{n^*}} (F_{lh(\eta)} )$ with the second player following
his winning strategy, it's easy to see that for each $\eta \in
T_{\aleph_{n^*}}$ we have $|\Suc_{T_{\aleph_{n^*}}} (\eta ) | = \aleph_{lh(\eta)}$.  Now let $b$ be an
infinite branch of $T_{\aleph_{n^*}}$ such that $b(i) >
\Sup(N_{\aleph_{n^*} ,
b\upharpoonright i} \cap \aleph_{m_i} )$, where $N_{\aleph_{n^*},
b\upharpoonright i}$ is defined in the
same way as $N_{\xi ,\eta}$ was defined above.  Such a branch $b$ clearly
exists.

Now, let $B = sk^A (\{ b(i) : i \in \omega \} \cap \{ \alpha_{\xi ,m} : \xi
< \aleph_{n^*} \wedge n^* < m \not\in A\})$.  Then for each $m \in A$, the
set $\la b(i) : i < \omega \wedge m_i = m\ra$ is cofinal in $B \cap
\aleph_m$.  Furthermore, for $n^* < m \not \in A$, $\langle \alpha_{\xi ,m}
: \xi < \aleph_{n*}\ra$ is cofinal in $B \cap \aleph_m$.  Hence $B$ is as
required.  \hfill   ${\square}$

\bigskip

\noindent {\bf Theorem 4.2.}  {\it Assume $GCH$.  Let $0 \leq m < n^* <
\omega$ and $E\subseteq \omega \backslash (n^* + 1)$ be such that for all $i
\in E, i + 1 \not \in E$ and let $j(i) = \Max (E\cap i)$.
Let $\la n_i : i \in E\ra$ with each $n_i < n^*$.
Suppose that for each $i \in E$, there is an $\aleph_i$-complete filter
$F_i$ on $\aleph_i$ containing all clubs of $\aleph_i$ such that $W_i = \{
\alpha < \aleph_i : cf(\alpha ) = \aleph_{n_i}\} \in F^+_i$ and the second
player has a winning strategy in the same $GM'_{\aleph_m , \aleph_i} (F_i).$

Let $f: \omega \rightarrow n^*$ be the function defined by
$$
f(i) = \cases {
   n_i & if $i \in E$,\cr
   m   & if $n^* < i \not\in E$\cr}
$$
Then the set $S' = \{ x \subset \aleph_\omega : |x| \leq \aleph_{n^*}$ and
$(\forall i > n^*) cf (x \cap \aleph_i ) = \aleph_{f(i)} \}$ is stationary
in $[\aleph_\omega ]^{<\aleph_{n^* + 1} } $.

}

\bigskip

\n {\bf Remarks.}~~~1.  Instead of $GCH$, it's enough to assume for
$i < j$ in
$E$, we have $2^{\aleph_i} < \aleph_j$.

2.  The assumption is consistent, but not so if we strengthen it using
$GM_{\aleph_m} (F_i)$.  (By \cite{Sh:542})

\bigskip

\n {\bf Proof.}  Let $\lambda > > \aleph_\omega$ be a regular cardinal and $A =
\langle \HH (\lambda ), \in, \triangleleft , < n_i\! : \! i \! \in \! E >,
(\tau_i)_{i <
\omega , \cdots} \ra$be a fully skolemized structure with skolem functions
closed under compositions, where $\triangleleft$ is a well-ordering on
$\HH(\lambda )$.  In order to prove the theorem, it suffices to show that
there exists $B \prec A$ such that $|B| = \aleph_{n^*}$, $cf(B \cap
\aleph_i) = \aleph_{n_i}$ for each $i \in E$ and $cf (B \cap \aleph_i) =
\aleph_m $ for $n^* < i
\not \in E$.

For each $i > n^*$, let $h_i : W_i \rightarrow [ \aleph_i ]^{<\aleph_{n^*}}$
be defined by $h_i (\delta ) = X_{i,\delta}$, where $X_{i,\delta}$ is the
$\triangleleft$-least cofinal subset of $\delta$ of cardinality
$\aleph_{n_i}$.  Note that each $h_i$ is definable in $A$.

We now define $\la A_{i,0} , B_{i,\xi} : i \in E, \xi < \aleph_m\ra$ and
$\la A_\xi : \xi \in \aleph_m \ra$ by induction on $\xi < \aleph_m$ such
that

\item {1)}  For each $i \in E, \la A_{i,0}, P_{i,\xi}, B_{i,\xi} : \xi <
\aleph_m \ra$ is a play of the game $GM'_{\aleph_m} (F_i)$ with the second
player following his winning strategy;

\item {2)}  $A_{i,0} = W_i, A_0 = sk^A  (\{ \emptyset \})$ and
$A_\xi = \cup_{\xi ' < \xi} A_{\xi'}$ if $\xi$ is limit;

\item {3)}  $A_\xi \prec A, |A_\xi | < \aleph_{n^*}$ and $A_\xi \subseteq
A_{\xi + 1}$;

\item {4)}  for each $n^* < j < \omega , j \not\in E$ we here
$\sup (A_\xi \cap \aleph_j) < \sup (A_{\xi + 1} \cap \aleph_j );$

\item {5)}  for each $i \in E$, for all
$\delta \in B_{i,\xi + 1}, $ we have $sk^A (A_\xi \cup X_{i,\delta  } ) \cap \aleph_{i-1}
\subseteq A_{\xi + 1}$ and $sk^A (A_\xi \cup \delta) \cap \aleph_i =
\delta$.

We simulate the games $GM'_{\aleph_m, \aleph_i} (F_i)$ for $i \in E$
simultaneously.
The first player chooses, $A_{i,0} = W_i$ and then divides it into less than
$\aleph_i$ parts for his (or her) $0^{th}$-move in the game $GM'_{\aleph_m,
\aleph_i} (F_i)$.  The second player always follows the winning strategy.
For successor stage, suppose we have constructed $\la B_{i,\xi} : i \in
E\ra$ and $A_\xi$.  For
$i \in E$, let
$C_i = \{ \alpha < \aleph_i : sk^A (\alpha \cup A_\xi) \cap \aleph_i =
\alpha \}$.  Then $C_i$ is a club in $\aleph_i$.  So $C_i \in F_i$.  For
each $i \in E$, consider the function $f_i : B_{i,\xi} \cap C_i \rightarrow
[\aleph_{i-1} ]^{<\aleph_{n^*}}$ defined by $f_i (\delta ) = sk^A (A_\xi \cup
X_{i,\delta} ) \cap \aleph_{i-1}$.  The first player divides $B_{i,\xi}$
into $\aleph_{i-1}$ parts as follows:
$P_{i,\xi} = \{ f^{-1}_i \{ x\} : x \in [\aleph_{i-1} ]^{<\aleph_{n^*} } \}
\cup \{ B_{i,\xi} \backslash C_i \}$.  (Note that $| [\aleph_{i-1}
]^{<\aleph_{n^*} } | = \aleph_{i-1}$ by GCH.)  The second player chooses one
of the parts, say $B_{i,\xi + 1}$ according to his winning strategy.  Note
that the second player will not choose $B_{i , \xi} \backslash C_i$ as his
move since $B_{i,\xi} \backslash C_i = \emptyset \mod(F_i)$.  (Otherwise he will
lose right away).  So there must be some $X_i \in [\alpha_{i-1}
]^{<\aleph_{n^*} }$ such that $f''_i B_{i,\xi + 1} = \{ X_i\}$.  Now let
$X = \cup_{i\in E} X_i$ and $\alpha_j = \sup (A_\xi \cap \aleph_j)$ for
$j\not \in E$ and  $A_{\xi + 1} = sk^A (A_\xi \cup X \cup \{\alpha_j : j
\not\in E\})$.

For limit stage, having defined $\la B_{i,\xi '} : \xi' < \xi\ra$, the first
player just divides $\cap_{\xi ' < \xi} B_{i,\xi'}$ into
$\aleph_{i-1}$ parts anyway he wants.  We let $B_{i,\xi}$ be the move of the
second player following his winning strategy.  This completes the
construction and the sequences $\la A_{i,0}, B_{i,\xi } : i \in E, \xi \in
\aleph_m \ra$ and $\la A_\xi : \xi \in \aleph_m\ra$ clearly satisfy clauses
(1)-(5) above.

Now, let $A^* = \cup_{\xi < \aleph_m} A_\xi$ and $W'_i = \cap_{\xi <
\alpha_m} B_{i,\xi}$.  Then $cf( A^* \cap \aleph_j) = \aleph_m$ for
$n^* < j \not\in E$ by clause (4) above, and each $W'_i $ is unbound in
$\aleph_i$.

Let's enumerate $E$ as $\la i_n : n \in \omega\ra$ in increasing order.  We
choose $\delta_{i_n} \in W'_{i_n}$.  Let $B_n = sk^A(A^* \cup \cup_{k\leq n}
X_{i_k,\delta_{i_k}}$ and $B'_n = sk^A (A^* \cup X_{i_n, \delta_{i_n}}).$
Then we have that $\sup (B'_n \cap \aleph_{i_n}) = \delta_{i_n}$ by clause
(5) above.  Also, we have that $B'_n \cap \aleph_{i_n-1} \subseteq A^*$
since if $\alpha \in B'_n \cap \aleph_{i_n-1}$ then $\alpha \in sk^A (A_\xi
\cup X_{i_n, \delta_{i_n}}) \cap \aleph_{i_n-1} \subseteq A_{\xi + 1}
\subseteq A^*$ for some $\xi < \aleph_m$.

\bigskip

\n {\bf Claim.}  {\it For all $n < \omega$, we have

\item {\rm a)}  $B_n \cap \aleph_{i_n-1} = B_{n-1} \cap \aleph_{i_n-1} $ for $n
> 0$

\item {\rm b)}  $\sup (B_n \cap \aleph_{i_n}) = \delta_{i_n} $

\item {\rm c)}  $( \forall i_0 < j \not\in E)
[\sup (B_n \cap \aleph_j) = \sup
(A^* \cap \aleph_j )]$

}

\bigskip

\n {\bf Proof.}  To prove a), it suffices to show that for any $\alpha \in B_n
\cap \aleph_{i_n-1}, \alpha \in B_{n-1}$.  Let $\alpha \in B_n \cap
\aleph_{i_n-1}$.  For simplicity, we may assume $\alpha = \tau(a^*, x_0,
\cdots , x_{n-1}, x_n)$ for $a^* \in A^*, x_k \in X_{i_k,\delta_{i_k}}$ for
$k \leq n$ and for some skolem function $\tau$.

Let $f: \prod_{k<n} \aleph_{i_k} \rightarrow
\aleph_{i_n -1}$ be the function defined by letting $f(\vec \beta ) = \tau
(a^*, \vec \beta, x_n)$ if $\tau (a^*, \vec \beta, x_n) < \aleph_{i_n-1}$
and $f(\vec B) = 0$ otherwise.  Then $f$ is definable from $a^*$ and $x_n$.
So $f\in B'_n.$

Now, let $\vec f = \la f_\xi : \xi < \aleph_{i_n-1} \ra$ be a list of all
the functions from $\prod_{k<n} \aleph_{i_k}$ to $\aleph_{i_n-1}$.  (Note
this is possible by GCH and $i_{n-1} < i_n -1.$)  By definability, we can
choose $\vec f \in A^*$.  But then $B'_n \models ( \exists \xi <
\aleph_{i_n-1} ) f_\xi = f$.  Let $\xi \in B'_n \cap \aleph_{i_n-1}$ be such
that $f_\xi = f$.  Then $\xi \in A^* \subseteq B_{n-1}$ since $B'_n \cap
\aleph_{i_n-1} \subseteq A^*$.  So $f = f_\xi \in B_{n-1}$.  Therefore,
$\alpha = f(x_0, \ldots , x_{n-1} ) \in B_{n-1} $ since $x_k \in B_{n-1}$
for all $k < n$.  We have thus proved part a) of the claim.

Clause b) follows from Lemma 2.1 and $\sup (B'_n \cap \aleph_{i_n}) =
\delta_{i_n}$.  By Lemma 2.1, $\sup (B_n \cap \aleph_{i_n} ) = \sup(B'_n
\cap \aleph_{i_n} ) = \delta_{i_n} .$)

We prove c) by induction on $n$.  If $n = 0,$ clause c) follows from Lemma 2.1.
Now suppose c) holds for $n - 1$.  We want to show c) holds for $n$.  By a)
and induction hypothesis, $\sup (B_n \cap \aleph_j) = \sup (B_{n-1} \cap
\aleph_j) = \sup (A^* \cap \aleph_j)$ if $i_0 < j < i_n$ and $j \not\in E$.
For $i_n < j \not\in E, \sup (B_n \cap \aleph_j ) = \sup (B'_n \cap
\aleph_j) = \sup (A^* \cap \aleph_j)$ by Lemma 2.1.  This finishes the proof
of the claim.   \hfill  ${\square}$

\bigskip

We now can complete the proof of Theorem 4.2.  Let $B^* = \cup_{n<\omega}
B_n$.  Then $B^* \prec A$ and $(\forall i \geq \min (E))[ cf (B^* \cap
\aleph_i) = \aleph_{f(i)} ] $ by the above claim.  (Note that if $i > n^*$ and
$i\leq i_n$, then $B^* \cap \aleph_i = B_n \cap \aleph_i$ by the claim.)

Finally, let $B = sk^A(B^* \cup \aleph_{n^*})$.  $B$ is as required again by
Lemma 2.1.  So we have finished the proof of Theorem 4.2.  \hfill
${\square}$

\bigskip

\noindent {\bf Concluding Remarks.}  The most natural context (at least for the second author) is having a
constant cardinal $\kappa$, set ${\frak a}$ of regular cardinals.  Let
$\lambda = \sup ({\frak a})$ and we look for stationary subsets of
$[\lambda]^{<\kappa}$.  Let $\FF^\kappa_{\frak a} = \{ f : f $ is a
function with domain ${\frak a}$ and $f (\theta )$ is a regular cardinal
$< \theta $ and $\sup \hbox {Rang} (f)  < \kappa \}$.  For an ideal $J$ on
${\frak a}$
and $f \in \FF^\kappa_{\frak a}$ we define $S^J_f
= \{ A \subseteq \lambda: |A| < \kappa$ and for some ${\frak b} \in J$ for
each $\theta \in {\frak a} \backslash {\frak b}$ the set $A \cap \theta$ is
a bounded subset of $\theta$ with order type of cofinality $f(\theta)\}$
.  Note that ${\frak a}, \lambda$ and $\kappa$
can be reconstructed from $f$ so we can just say ``$S_f = S^J_f$ is
stationary''.  We
call the framework simple if $\kappa \leq \Min({\frak a})$, and we
concentrate on it.  If $J = \{ \emptyset \}$ we may omit it.

\item {$(*)_1$}  If ${\frak a}$ has a maximal element $\theta $ and $f \in
\FF^\kappa_{\frak a}$  then ($f\upharpoonright ({\frak a} \backslash \{ \theta \})
\in \FF_{  {\frak a} \backslash \{ \theta \}}$ and) $S^J_f$ is stationary iff
$S_{f\upharpoonright({\frak a} \backslash \{ \theta \})}$ is stationary.

\medskip

\item {$(*)_2$}  for $f \in \FF^\kappa_{\frak a}$ and $\theta$ we have
$S^J_f$ is stationary iff $S^J_{f\upharpoonright ({\frak a} \cap \theta)}$
is stationary and $S^J_{f\upharpoonright ({\frak a} \backslash \theta )}$ is
stationary.

\medskip

\item {$(*)_3$} assume ${\frak a}$ has no last element and $f \in
\FF^\kappa_{\frak a}$, then $S_f$ is stationary iff $(a) + (b)$ where:

$(a)$ for
every algebra $M$ with universe $\sup ({\frak a})$ for some $N \prec M$ we
have:  for every $\theta \in {\frak a}$ large enough,
$$
cf(\sup (N \cap \theta ) = f(\theta ) .
$$

$(b)$ for every $\theta \in {\frak a}$ the set $S^J_{f\upharpoonright
({\frak a} \cap \theta )}$ is stationary.

\medskip

\item {$(*)_4$}  Assume $n^* < \omega$ and $\lambda_1 <
\lambda_2 < \cdots < \lambda_{n^*}$ are
member of $pcf ({\frak a})$ which are $> \sup ({\frak a})$.  Assume further
$f \in \FF^\kappa_{\frak a}$, $\la {\frak b}_1 , \cdots , {\frak b}_{n^*}
\ra$ is a partition of ${\frak a}$, and $f\upharpoonright {\frak b}_{e}$ is
constant and $\lambda_e = tcf (\prod {\frak b}_e, <_{J\upharpoonright {\frak
b}_e} )$ $\underline {\rm then}$ $S^J_f$ is stationary.

[Why?  by the proof of Theorem 3.1.  I.e. by $(*)_2$ used several times
${\rm wlog~} \Min ( {\frak a}) > | {\frak a} |^{n^*+2}$, then by the proof of 3.1
if $1 \leq i \leq n^* \Rightarrow f(i) > |{\frak a} |$ we succeed.  Lastly use $(*)_3$
possibly several times.]

\medskip

\item {$(*)_5$}  Assume (a)  ${\frak b} \subseteq {\frak a}$ is countable,
$f \in \FF^\kappa_{\frak a}$,  $f \upharpoonright {\frak b}$ is constantly
$\aleph_0$ and $f\upharpoonright ({\frak a} \backslash {\frak b})$ is
constantly $\sigma$, $\delta$ is a unit ordinal of cofinably $\sigma$ but
$< \kappa$ and $\delta$ is divisible by $| {\frak a} |$.

\item {}(b)  for $\theta \in {\frak a}$, the second player has a winner
strategy in the game $GM_{\sigma_0} (F_\theta , \theta \cap ({\frak a}
\backslash {\frak b}))$ (see second remark to Theorem 4.1).

$\underline {\rm Then}$
$S_f$  is stationary.

\item {}[Why?  repeat the proof of Theorem 4.1, but we let $\la \theta_i : i < \omega
\ra$ list ${\frak b}$, each appearing infinitely often, and $\TT = \{ T: T $
subtree of $\cup_{e< \omega } \prod_{i<e} \theta_i$ such that for every
$\eta\in T$ of length $i$ we have
$\{ \alpha < \theta_i :  \eta^{\hat{} }  \la i\ra \in
T\} \not = \emptyset \mod F_{\theta_i} \}$ let $\la (\theta_\xi, m_\xi) : \xi
< \delta \ra$ be such that: $\theta_\alpha \in {\frak a} \backslash {\frak b},
m_\alpha < \omega$, and each such pair occurs boundedly often.  Then define
the $T_\xi \in \TT $ as in the proof of Theorem 4.1, in $T_{\xi + 1}$ we take
care of every $\eta \in T_{\xi + 1}$ of length $\leq m_\xi$.]

\medskip

\item {$(*)_6$}  Assume (a)  ${\frak b} \subseteq {\frak a}$, $\sigma = cf
(\sigma ) < \kappa, f\in \FF^\kappa_{\frak a}$, $f\upharpoonright ({\frak
a} \backslash {\frak b})$ is constantly $\sigma$, $\delta$ is an ordinal $<
\kappa$ of cofinality $< \sigma $ and let $\sigma_\theta = [(\sup ({\frak a} \cap
\theta ))^{\sigma + \sup ({\frak b} \cap \theta )} ]^+$.

\item ~~~~~{(b)}  for $\theta \in {\frak b}$, $F_\theta$ is a $\sigma_\theta$-complete
filter on $\theta$ extending the club filter such that player two has a winning
strategy in the game $GM'_{\delta,\sigma_\theta} (F_\theta).$

$\underline
{\rm Then}$ $S_f$ is stationary.

\item {}[Why?  by the proof of Theorem 4.2]

\vfill\eject

REFERENCES.  
\bigskip

\bibliographystyle{lit-plain}
\bibliography{lista,listb,listx}

\shlhetal

\bye